\documentclass[a4paper]{jpconf}
\usepackage{amsfonts,amsthm,amsmath,amssymb,
	amsbsy,euscript}

\newtheorem{theor}{Theorem}

\theoremstyle{definition}

\newtheorem{lemma}[theor]{Lemma}
\newtheorem{define}{Definition}

\newtheorem*{notation}{Notation}
\newtheorem{problem}{Problem}
\newtheorem{open}[problem]{Open problem}
\newtheorem{example}{Example}
\theoremstyle{remark}
\newtheorem{rem}{Remark}

\newcommand{\pinner}{\mathbin{\mathchoice
{\hbox{\vrule width0.6em depth0pt height0.4pt
	\vrule width0.4pt depth0pt height0.8ex}}
{\hbox{\vrule width0.6em depth0pt height0.4pt
	\vrule width0.4pt depth0pt height0.8ex}}
{\hbox{\kern0.14em
	\vrule width0.48em depth0pt height0.4pt
	\vrule width0.4pt depth0pt height0.6ex\kern0.14em}}
{\hbox{\kern0.1em
	\vrule width0.39em depth0pt height0.4pt
	\vrule width0.4pt depth0pt height0.5ex\kern0.1em}}}}

\newcommand{\cN}{\mathcal{N}}

\newcommand{\cP}{\mathcal{P}}
\newcommand{\cQ}{\mathcal{Q}}

\newcommand{\cX}{{\EuScript X}}    

\newcommand{\bx}{{\boldsymbol{x}}}

\newcommand{\veps}{\varepsilon}

\newcommand{\Id}{{\mathrm d}}

\DeclareMathOperator{\Alt}{Alt}

\DeclareMathOperator{\Jac}{Jac}

\newcommand{\schouten}[1]{\lshad {#1} \rshad}

\newcommand{\by}[1]{\textit{{#1}}}
\newcommand{\jour}[1]{\textit{{#1}}}
\newcommand{\vol}[1]{\textbf{{#1}}}
\newcommand{\book}[1]{\textrm{{#1}}}

\newcommand*{\vcenteredhbox}[1]{\begingroup
\setbox0=\hbox{#1}\parbox{\wd0}{\box0}\endgroup}

\DeclareSymbolFont{extraup}{U}{zavm}{m}{n}
\DeclareMathSymbol{\varheartsuit}{\mathalpha}{extraup}{86}
\DeclareMathSymbol{\vardiamondsuit}{\mathalpha}{extraup}{87}

\def\oldvec{\mathaccent "017E\relax } 
\DeclareMathOperator{\Ori}{\mathsf{O\oldvec{r}}}

\begin{document}\pagestyle{plain}

\title{Infinitesimal deformations of Poisson bi\/-\/vectors\\ 
using the Kontsevich graph calculus}

\author{Ricardo Buring${}^1$, Arthemy V Kiselev${}^2$ and Nina Rutten${}^2$}

\address{${}^1$ 
Institut f\"ur Mathematik, 
Johannes Gutenberg\/--\/Uni\-ver\-si\-t\"at,
Staudingerweg~9, 
\mbox{D-\/55128} Mainz, Germany}
\address{${}^2$ Johann Ber\-nou\-lli Institute for Mathematics and Computer Science, University of Groningen,
P.O.~Box 407, 9700~AK Groningen, The Netherlands}

\ead{rburing@uni-mainz.de, A.V.Kiselev@rug.nl}

\begin{abstract}
Let $\cP$~be a Poisson structure on a 
finite\/-\/dimensional affine real manifold. Can $\cP$~be deformed in such a way that it stays Poisson\,? The language of Kontsevich graphs provides a universal approach --\,with respect to all affine Poisson manifolds\,-- to finding a class of solutions to this deformation problem. For that reasoning, several types of 
graphs are needed. In this paper we outline the algorithms to generate 
those graphs.
The graphs that encode deformations are classified by the number of internal vertices~$k$; 
for~$k \leqslant 4$ we present all solutions of the deformation problem.
For $k\geqslant5$, first reproducing the pentagon\/-\/wheel picture 
suggested at~$k=6$ by Kontsevich and Willwacher, we 
construct the heptagon\/-\/wheel cocycle that yields a new unique solution without $2$-\/loops and tadpoles at~$k=8$.
\end{abstract}

\pagestyle{plain}

\paragraph{\textbf{\textup{Introduction.}}}
This paper contains a set of algorithms to generate the Kontsevich graphs 
that encode polydifferential operators --\,in particular multi\/-\/vectors\,-- on Poisson manifolds.
We report a result of implementing such algorithms in the problem of 
fin\-ding symmetries 
of Poisson structures.
Namely, continuing the line of reasoning from 
\cite{tetra16,f16}, we find all the solutions of this deformation problem that are expressed by the Kontsevich graphs with at most four internal vertices. Next, we present 
one six\/-\/vertex solution (based on the previous work by Kontsevich~\cite{Kontsevich2017private} 
and 
Willwacher~\cite{Willwacher2017private}).
Finally, we find a heptagon\/-\/wheel eight\/-\/vertex graph 
which, after the orientation of its edges, gives 
a new universal Kontsevich flow.
We refer to~\cite{Ascona96, Kontsevich2017Bourbaki} for motivations, to~\cite{f16,cpp} for an exposition
of basic theory, and to~\cite{WeFactorize5Wheel} and~\cite{JNMP2017} for more details about the pentagon\/-\/wheel ($5+1$)-\/vertex and heptagon\/-\/wheel ($7+1$)-\/vertex
solutions respectively.
Let us remark that all the algorithms outlined here can be used without 
modification in the course of constructing all~$k$-\/vertex Kontsevich graph solutions 
with higher~$k \geqslant 5$ in the deformation problem under study.

\paragraph{\textbf{\textup{Basic concept.}}}
We work with real vector spaces 
generated by finite graphs of the following two types: 
(1)~$k$-\/vertex non-oriented graphs, without multiple edges nor tadpoles, endowed with a wedge ordering of edges, e.g., $E=e_1\wedge\dots\wedge e_{2k-2}$; 
(2)~oriented graphs on $k$ internal vertices and $n$ sinks such that every internal vertex is a tail of two edges with a given ordering Left $\prec$ Right. Every connected component of a non-oriented graph $\gamma$ is fully encoded by an ordering~$E$ on the set of adjacency relations for its vertices.\footnote{\label{FootZeroGraph}%
The edges are antipermutable so that a graph which equals minus itself 
--\,under a symmetry that induces a parity\/-\/odd permutation of edges\,--
is proclaimed to be equal to zero.
In particular (view $\bullet\!\!\!-\!\!\!\bullet\!\!\!-\!\!\!\bullet$), every graph possessing a symmetry which swaps an odd number of edge pairs is a zero graph. For example, the $4$-\/wheel $\mathsf{12}\wedge \mathsf{13}\wedge \mathsf{14}\wedge \mathsf{15}\wedge \mathsf{23}\wedge \mathsf{25}\wedge \mathsf{34}\wedge \mathsf{45} = I\wedge\cdots\wedge VIII$ or the $2\ell$-\/wheel at any $\ell>1$ is such; here, the 
reflection symmetry is $I\rightleftarrows III$,\ $V\rightleftarrows VII$, 
and~$VI \rightleftarrows VIII$.}
Every such oriented graph is given by the list of ordered pairs of directed edges. An edge swap $e_i\wedge e_j = -e_j\wedge e_i$ and the reversal 
Left $\leftrightarrows$ Right of those edges' order in the tail vertex implies the change of sign in front of the graph at hand.\footnote{An oriented graph equals minus itself, hence it is a zero graph if there is a permutation of labels for its internal vertices such that the adjacency tables for the two vertex labellings coincide but the two realisations of the same graph differ by the ordering of outgoing edges at an odd number of internal vertices (see Example~\ref{ExZeroGraph} below).} 

\begin{example}\label{Ex5Wheel}
The sum $\boldsymbol{\gamma}_5$ of two 6-vertex 10-edge graphs,
\begin{align*}
\boldsymbol{\gamma}_5 = \mathsf{12}^{(I)}\wedge \mathsf{23}^{(II)}\wedge \mathsf{34}^{(III)}\wedge \mathsf{45}^{(IV)} \wedge \mathsf{51}^{(V)} \wedge \mathsf{16}^{(VI)} \wedge \mathsf{26}^{(VII)}\wedge \mathsf{36}^{(VIII)}\wedge \mathsf{46}^{(IX)}\wedge \mathsf{56}^{(X)}&\\
 + \tfrac{5}{2}\cdot \mathsf{12}^{(I)}\wedge \mathsf{23}^{(II)}\wedge \mathsf{34}^{(III)} \wedge \mathsf{41}^{(IV)}\wedge \mathsf{45}^{(V)}\wedge \mathsf{15}^{(VI)}\wedge \mathsf{56}^{(VII)}\wedge \mathsf{36}^{(VIII)}\wedge \mathsf{26}^{(IX)}\wedge \mathsf{13}^{(X)}&,
\end{align*}
is drawn in Theorem~\ref{ThMainMany} on p.~\pageref{Where5WheelIs} 
below.
\end{example}

\begin{example}\label{Ex3Wheel}
The sum $\cQ_{1:\frac{6}{2}}$ of three oriented 8-edge graphs on $k=4$ internal vertices and $n=2$ sinks (enumerated using 0 and 1, see the notation on p.~\pageref{DefEncodingKgraph}),
$$\cQ_{1:\frac{6}{2}} = \mathsf{2\ 4\ 1\ \ \ 0\ 1\ \ 2\ 4\ \ 2\ 5\ \ 2\ 3} - 3(\mathsf{2\ 4\ 1\ \ \ 0\ 3\ \ 1\ 4\ \ 2\ 5\ \ 2\ 3} + \mathsf{2\ 4\ 1\ \ \ 0\ 3\ \ 4\ 5\ \ 1\ 2\ \ 2\ 4})$$
is obtained from the non-oriented tetrahedron graph $\boldsymbol{\gamma}_3 = 
\mathsf{12}^{(I)}\wedge \mathsf{13}^{(II)}\wedge \mathsf{14}^{(III)}\wedge \mathsf{23}^{(IV)}\wedge \mathsf{24}^{(V)}\wedge \mathsf{34}^{(VI)}$ on four vertices and six edges by taking all the admissible edge orientations 
(see Theorem~\ref{ThOrient} and Remark~\ref{RemSigns}
). 
\end{example}

\paragraph{\textbf{\textup{I.1.}}}
Let $\gamma_1$ and $\gamma_2$ be connected non-oriented graphs. The definition of insertion $\gamma_1\circ_i\gamma_2$ 
of the entire graph $\gamma_1$ into vertices of $\gamma_2$ and the construction of Lie bracket $[\cdot,\cdot]$ of graphs and differential~$\Id$ in the non\/-\/oriented graph complex, referring to a sign convention, are as follows
(cf.~\cite{Ascona96} and~\cite{DolgushevRogersWillwacher,KhoroshkinWillwacherZivkovic,WillwacherGRT}); these definitions apply to sums of graphs by linearity.

\begin{define}\label{DefInsert}
The insertion $\gamma_1\circ_i\gamma_2$ of a $k_1$-vertex graph $\gamma_1$ with ordered set of edges $E(\gamma_1)$ into a graph $\gamma_2$ with $\# E(\gamma_2)$ edges on $k_2$ vertices is a 
sum of 
graphs on $k_1+k_2-1$ vertices and $\# E(\gamma_1)+\# E(\gamma_2)$ edges. 
Topologically, the sum $\gamma_1\circ_i\gamma_2 = \sum(\gamma_1\rightarrow v \text{ in }\gamma_2)$ 
 consists of all the graphs in which a vertex $v$ from $\gamma_2$ is replaced by the entire graph $\gamma_1$ and the edges tou\-ch\-ing $v$ in $\gamma_2$ are re-attached to the vertices of $\gamma_1$ in all possible ways.\footnote{Let the enumeration of vertices in every such term in the sum start running over the enumerated vertices in $\gamma_2$ until $v$ is reached.
Now the enumeration counts the vertices in the graph $\gamma_1$ and then it resumes with the remaining vertices (if any) that go after~$v$ 
in~$\gamma_2$.}
By convention, in every new term the edge ordering is~$E(\gamma_1)\wedge E(\gamma_2)$. \end{define}

To simplify sums of graphs, first eliminate the zero graphs.
Now suppose that in a sum, two non\/-\/oriented graphs, say $\alpha$ and $\beta$, are isomorphic (topologically, i.e.\ regardless of the respective vertex labellings and edge orderings $E(\alpha)$ and~$E(\beta)$).
By using that isomorphism, which establishes a 1--1~correspondence between the edges, extract the sign from the equation $E(\alpha) = \pm E(\beta)$.
If~``+'', then $\alpha = \beta$; else $\alpha = -\beta$.
Collecting similar terms is now elementary.

\begin{lemma}
The bi\/-\/linear graded skew\/-\/symmetric operation,
\[
[\gamma_1,\gamma_2] = \gamma_1\circ_i\gamma_2 - (-)^{\# E(\gamma_1)\cdot\# E(\gamma_2)} \gamma_2\circ_i\gamma_1,
\]
is a Lie bracket on the vector space $\mathfrak{G}$ of non-oriented graphs.\footnote{The postulated precedence or antecedence of the wedge product of edges from $\gamma_1$ with respect to the edges from $\gamma_2$ in every graph within $\gamma_1\circ_i\gamma_2$ produce the operations $\circ_i$ which coincide with or, respectively, differ from Definition~\ref{DefInsert} by the sign factor $(-)^{\# E(\gamma_1)\cdot\# E(\gamma_2)}$. The same applies to the Lie bracket of graphs $[\gamma_1,\gamma_2]$ if the operation $\gamma_1\circ_i\gamma_2$ is 
the insertion of $\gamma_2$ into $\gamma_1$ (as in~\cite{KhoroshkinWillwacherZivkovic}). Anyway, the notion of $\Id$-cocycles which we presently recall is well defined and insensitive to such sign~ambiguity.}
\end{lemma}

\begin{lemma}
The operator $\Id(\text{graph}) = [\bullet\!\!\!-\!\!\!\bullet, \text{graph}]$ is a differential: $\Id^2 = 0$.
\end{lemma}

In effect, the mapping $\Id$ blows up every vertex $v$ in its argument in such a way that whenever the number of adjacent vertices $N(v) \geqslant 2$ is sufficient, each end of the inserted edge $\bullet\!\!\!-\!\!\!\bullet$ is connected with the rest of the graph by at least one edge.

Summarising,
the real vector space~$\mathfrak{G}$ of non\/-\/oriented graphs is a differential graded Lie algebra \textup{(}dgLa\textup{)} with 
Lie bracket $[\cdot,\cdot]$ and differential~$\Id = [\bullet\!\!\!-\!\!\!\bullet, \cdot]$.
The graphs $\boldsymbol{\gamma}_5$ and $\boldsymbol{\gamma}_3$ from Examples~\ref{Ex5Wheel} and~\ref{Ex3Wheel} are $\Id$-\/cocycles.
Neither is exact, hence marking a nontrivial cohomology class in the non-oriented graph complex.

\begin{theor}[{{\cite[Th.\:5.5]{DolgushevRogersWillwacher}}}]\label{ThWheelCocycles}
At every $\ell \in \mathbb{N}$ in the connected graph complex there is a $\Id$-cocycle on $2\ell + 1$ vertices and $4\ell+2$ edges.
Such cocycle contains the $(2\ell+1)$-wheel in which, by definition, the axis vertex is connected with every other vertex by a spoke so that each of those $2\ell$ vertices is adjacent to the axis and two neighbours\textup{;} the cocycle marked by the $(2\ell+1)$-wheel graph can contain other $(2\ell+1, 4\ell+2)$-\/graphs \textup{(}see Example~\textup{\ref{Ex5Wheel}} and~\textup{\cite{JNMP2017})}.
\end{theor}

\paragraph{\textbf{\textup{I.2.}}} 
The \emph{oriented} graphs under study are built over $n$ sinks from $k$ wedges $\xleftarrow{i_\alpha} \bullet \xrightarrow{j_\alpha}$ (here $\xleftarrow{i_\alpha} \prec \xrightarrow{j_\alpha}$) so that every edge is decorated with its own summation index which runs from $1$ to the dimension of a given affine Poisson manifold $(\cN, \cP)$.
Each edge $\xrightarrow{i}$ encodes the derivation $\partial/\partial x^i$ of the arrowhead object with respect to a local coordinate $x^i$ on $\cN$.
By placing an $\alpha$th copy $P^{i_\alpha j_\alpha}(\bx)$ of the Poisson bi-vector $\cP$ in the wedge top $(1 \leqslant \alpha \leqslant k$), by taking the product of contents of the $n+k$ vertices (and evaluating all objects at a point $\bx \in \cN$), and summing over all indices, we realise a polydifferential operator in $n$ arguments; the operator coefficients are differential-polynomial in~$\cP$.
Totally skew-symmetric operators of differential order one in each argument are well-defined $n$-vectors on the affine manifolds $\cN$.

The space of multi\/-\/vectors $G$ encoded by oriented graphs is equipped with a graded Lie algebra structure, namely the Schouten bracket $\schouten{\cdot,\cdot}$.
Its realisation in terms of oriented graphs is shown in \cite[Remark 4]{f16}.
Recall that by definition the bi-vectors $\cP$ at hand are Poisson by satisfying the Jacobi identity $\schouten{\cP, \cP} = 0$.
The Poisson differential $\partial_\cP = \schouten{\cP, \cdot}$ now endows the space of multi-vectors on $\cN$ with the differential graded Lie algebra (dgLa) structure.
The cohomology groups produced by the two dgLa structures introduced 
so far are correlated 
by the edge orientation mapping~$\Ori$.

\begin{theor}[{{\cite{Ascona96} and \cite[App.\:K]{WillwacherGRT}}}]%
\label{ThOrient}
Let $\gamma \in \ker \Id$ be a cocycle on $k$~vertices and $2k-2$~edges in the non-oriented graph complex. Denote by~$\{\Gamma\} \subset G$ the subspace spanned by all those bi\/-\/vector graphs~$\Gamma$ which are obtained from \textup{(}each connected component in\textup{)} $\gamma$ by adding to it two edges to the new sink vertices and then by taking the sum of graphs with all the admissible orientations of the old $2k-2$~edges
\textup{(}so that a set of Kontsevich graphs built of 
$k$~wedges is produced\textup{)}.
Then in that subspace $\{\Gamma\}$ there is a sum of graphs that encodes a nonzero Poisson cocycle $Q(\cP) \in \ker \partial_\cP$.
\end{theor}


Consequently, to find some cocycle $Q(\cP) 
$ in the Poisson complex on any affine Poisson manifold it suffices to find a cocycle in the non-oriented graph complex and then consider 
the sum of graphs which are produced by the orientation mapping~$\Ori$. 
On the other hand, to list all the $\partial_\cP$-cocycles $Q(\cP)$ encoded by the bi-vector graphs made of 
$k$ wedges ${\leftarrow}{ \bullet }{\rightarrow}$, one must generate all the relevant oriented graphs and solve the equation $\partial_\cP(Q) \doteq 0$ via $\schouten{\cP,\cP} = 0$, that is, solve graphically the factorisation problem $\schouten{\cP, Q(\cP)} = \Diamond(\cP, \schouten{\cP,\cP})$ in which the cocycle condition in the left-hand side holds by virtue of the Jacobi identity in the right.
Such construction of some and classification (at a fixed $k > 0$) of all universal infinitesimal symmetries of Poisson 
brackets are the problems which we explore in this paper.

\begin{rem}\label{RemSigns}
To the best of our knowledge \cite{Kontsevich2017private}, in a bi-vector graph $Q(\cP) = \Ori(\gamma)$, at every internal vertex which is the tail of two oriented edges towards other internal vertices, the edge ordering Left $\prec$ Right 
is inherited from a chosen wedge product~$E(\gamma)$ of edges in the non-oriented 
graph~$\gamma$.
How are the new edges towards the sinks ordered, either between themselves at a vertex or with respect to two other oriented edges, coming from~$\gamma$ and issued from different vertices in~$Q(\cP)$\,?
Our findings in~\cite{WeFactorize5Wheel} will help us to verify the order preservation claim and assess answers to this question.
\end{rem}




\section{The Kontsevich graph calculus}\label{SecGraphs}

\begin{define}
\label{DefKgraph}
Let us consider a class of oriented graphs on $n+k$ vertices labelled $0$,\ $\ldots$,\ $n+k-1$ such that the consecutively ordered vertices $0$,\ $\ldots$,\ $n-1$ are sinks, and each of the internal vertices $n$,\ $\ldots$,\ $n+k-1$ is a source for two edges. 
For every internal vertex, the two outgoing edges are ordered using $L \prec R$: the preceding edge is labelled $L$ (Left) and the other is $R$ (Right).
An oriented graph on $n$ sinks and $k$ internal vertices 
is a \emph{Kontsevich graph} of type~$(n,k)$.
\end{define}

For the purpose of defining a graph normal form, we now consider a Kontsevich graph $\Gamma$ together with a {\em sign} $s \in \{0, \pm 1\}$, denoted by concatenation of the symbols: $s\Gamma$.

\begin{notation}[Encoding of the Kontsevich graphs]\label{DefEncodingKgraph}
The format to store a signed graph $s\Gamma$ for a Kontsevich graph $\Gamma$ is the integer number $n>0$, the integer $k \geqslant 0$, the sign $s$, followed by the (possibly empty, when $k=0$) list of $k$~ordered pairs of targets for edges issued from the internal vertices $n$,\ $\ldots$,\ $n+k-1$, respectively.
The full format is then ($n$,\ $k$,\ $s$; list of ordered pairs).
\end{notation}

\begin{define}[Normal form of a Kontsevich graph]
The list of targets in the encoding of a graph $\Gamma$ can be considered as a $2k$-digit integer written in base-$(n+k)$ notation.
By running over the entire group $S_k \times (\mathbb{Z}_2)^k$, and by this over all the different re\/-\/labellings of~$\Gamma$, we obtain many different integers written in base-$(n+k)$.
The {\em absolute value} $|\Gamma|$ of $\Gamma$ is the re\/-\/labelling of $\Gamma$ such that its list of targets is {\em minimal} as a nonnegative base-$(n+k)$ integer.
For a signed graph $s\Gamma$, the {\em normal form} is the signed graph $t|\Gamma|$ which represents the same polydifferential operator as $s\Gamma$.
Here we let $t=0$ if the graph is zero (see Example \ref{ExZeroGraph} below).
\end{define}

\begin{example}[Zero Kontsevich graph]\label{ExZeroGraph}
{\unitlength=1mm
\special{em:linewidth 0.4pt}
\linethickness{0.4pt}
\begin{picture}(0,0)(-60,19)    
\put(15.00,25.00){\vector(1,0){10.00}}
\put(25.00,25.00){\vector(-1,-3){5.00}}
\put(25.00,25.00){\vector(1,-3){5.00}}
\put(25.00,15.00){\vector(1,-1){4.67}}
\put(25.00,15.00){\vector(-1,-1){4.67}}
\put(15.00,25.00){\line(1,-1){6.67}}
\put(23,17){\vector(1,-1){2}}
\bezier{16}(23.00,17.00)(24.00,16.00)(25.00,15.00)
\put(15.00,25.00){\circle*{1}}
\put(25.00,25.00){\circle*{1}}
\put(25.00,15.00){\circle*{1}}
\put(30.33,9.67){\circle*{1}}
\put(19.67,9.67){\circle*{1}}
\put(13.33,25.00){\makebox(0,0)[rc]{4}}
\put(26.67,25.00){\makebox(0,0)[lc]{3}}
\put(20.00,25.67){\makebox(0,0)[cb]{\tiny$R$}}
\put(17.67,21.00){\makebox(0,0)[ct]{\tiny$L$}}
\put(25.00,13.33){\makebox(0,0)[ct]{2}}
\put(20.00,8){\makebox(0,0)[ct]{0}}
\put(30.00,7.67){\makebox(0,0)[ct]{1}}
\end{picture}
}
Consider the graph with the encoding\\[1pt]
\parbox{131.5mm}{{\tt 2 3 1\ \ 0 1 0 1 2 3}.
The swap of vertices $2 \rightleftarrows 3$ is a symmetry of this graph, yet it also swaps the ordered edges $(4\to 2) \prec (4 \to 3)$, producing a minus sign.
Equal to minus itself, this Kontsevich graph is zero.}
\end{example}

\begin{notation}
Every Kontsevich graph~$\Gamma$ on $n$~sinks (or every sum~$\Gamma$ of such graphs) yields the sum~$\Alt\Gamma$ of Kontsevich graphs which is totally skew\/-\/symmetric with respect to the $n$~sinks content~$s_1$, $\ldots$, $s_n$. Indeed, let
\begin{equation}\label{EqMakeSkew}
(\Alt \Gamma)(s_1,\ldots,s_n) = \sum\nolimits_{\sigma\in\mathcal{S}_n} 
(-)^\sigma\, \Gamma(s_{\sigma(1)},\ldots,s_{\sigma(n)}).
\end{equation}
Due to skew\/-\/symmetrisation, the sum 
of graphs~$\Alt\Gamma$ can contain zero graphs or~\mbox{repetitions.}
\end{notation}

\begin{example}[The Jacobiator]
The left-hand side of the Jacobi identity is a skew sum of Kontsevich graphs (e.g. it is obtained by skew-symmetrizing the first term)
\begin{equation}\label{EqJacFig}
\vcenteredhbox{
\raisebox{3.3mm}
[6.5mm][3.5mm]{
\unitlength=1mm
\special{em:linewidth 0.4pt}
\linethickness{0.4pt}
\begin{picture}(12,15)
\put(0,-10){
\begin{picture}(12.00,15.00)
\put(0.00,10.00){\framebox(12.00,5.00)[cc]{$\bullet\ \bullet$}}
\put(2.00,10.00){\vector(-1,-3){1.33}}
\put(6.00,10.00){\vector(0,-1){4.00}}
\put(10.00,10.00){\vector(1,-3){1.33}}
\put(0.00,4.00){\makebox(0,0)[cb]{\tiny\it1}}
\put(6.00,4.00){\makebox(0,0)[cb]{\tiny\it2}}
\put(11.67,4.00){\makebox(0,0)[cb]{\tiny\it3}}
\end{picture}
}\end{picture}}}\ \ \ 
\mathrel{{:}{=}}
\text{\raisebox{-12pt}[25pt]{
\unitlength=0.70mm
\linethickness{0.4pt}
\begin{picture}(26.00,16.33)
\put(0.00,5.00){\line(1,0){26.00}}
\put(2.00,5.00){\circle*{1.33}}
\put(13.00,5.00){\circle*{1.33}}
\put(24.00,5.00){\circle*{1.33}}
\put(2.00,1.33){\makebox(0,0)[cc]{\tiny\it1}}
\put(13.00,1.33){\makebox(0,0)[cc]{\tiny\it2}}
\put(24.00,1.33){\makebox(0,0)[cc]{\tiny\it3}}
\put(7.33,11.33){\circle*{1.33}}
\put(7.33,11.33){\vector(1,-1){5.5}}
\put(7.33,11.33){\vector(-1,-1){5.5}}
\put(13,17){\circle*{1.33}}
\put(13,17){\vector(1,-1){11.2}}
\put(13,17){\vector(-1,-1){5.1}}
\put(3.00,10.00){\makebox(0,0)[cc]{\tiny$i$}}
\put(12.00,10.00){\makebox(0,0)[cc]{\tiny$j$}}
\put(24.00,10.00){\makebox(0,0)[cc]{\tiny$k$}}
\end{picture}
}}
{-}
\text{\raisebox{-12pt}[25pt]{
\unitlength=0.70mm
\linethickness{0.4pt}
\begin{picture}(26.00,16.33)
\put(0.00,5.00){\line(1,0){26.00}}
\put(2.00,5.00){\circle*{1.33}}
\put(13.00,5.00){\circle*{1.33}}
\put(24.00,5.00){\circle*{1.33}}
\put(2.00,1.33){\makebox(0,0)[cc]{\tiny\it1}}
\put(13.00,1.33){\makebox(0,0)[cc]{\tiny\it2}}
\put(24.00,1.33){\makebox(0,0)[cc]{\tiny\it3}}
\put(13,11.33){\circle*{1.33}}
\put(13,11.33){\vector(2,-1){10.8}}
\put(13,11.33){\vector(-2,-1){10.8}}
\put(18.5,17){\circle*{1.33}}
\put(18.5,17){\vector(-1,-1){5.2}}
\put(18.5,17){\vector(-1,-2){5.6}}
\put(13,15){\tiny $L$}
\put(17,12){\tiny $R$}
\put(4.00,10.00){\makebox(0,0)[cc]{\tiny$i$}}
\put(11.00,8.00){\makebox(0,0)[cc]{\tiny$j$}}
\put(22.00,10.00){\makebox(0,0)[cc]{\tiny$k$}}
\end{picture}
}}
{-}
\text{\raisebox{-12pt}[25pt]{
\unitlength=0.70mm
\linethickness{0.4pt}
\begin{picture}(26.00,16.33)
\put(0.00,5.00){\line(1,0){26.00}}
\put(2.00,5.00){\circle*{1.33}}
\put(13.00,5.00){\circle*{1.33}}
\put(24.00,5.00){\circle*{1.33}}
\put(2.00,1.33){\makebox(0,0)[cc]{\tiny\it1}}
\put(13.00,1.33){\makebox(0,0)[cc]{\tiny\it2}}
\put(24.00,1.33){\makebox(0,0)[cc]{\tiny\it3}}
\put(18.33,11.33){\circle*{1.33}}
\put(18.33,11.33){\vector(1,-1){5.5}}
\put(18.33,11.33){\vector(-1,-1){5.5}}
\put(13,17){\circle*{1.33}}
\put(13,17){\vector(-1,-1){11.2}}
\put(13,17){\vector(1,-1){5.1}}
\put(3.00,10.00){\makebox(0,0)[cc]{\tiny$i$}}
\put(13.00,10.00){\makebox(0,0)[cc]{\tiny$j$}}
\put(24.00,10.00){\makebox(0,0)[cc]{\tiny$k$}}
\end{picture}
}}
.
\end{equation}
The default ordering of edges is the one which we see. 
\end{example}

\begin{define}[Leibniz graph]
A \emph{Leibniz graph} is a graph whose vertices are either sinks, or the sources for two arrows, or the Jacobiator (which is a source for three arrows).
There must be at least one Jacobiator vertex.
The three arrows originating from a Jacobiator vertex must land on three distinct vertices.
Each edge falling on a Jacobiator works by the Leibniz rule on the two internal vertices in it.
\end{define}

\begin{example}
The Jacobiator itself is a Leibniz graph (on one tri-valent internal vertex).
\end{example}

\begin{define}[Normal form of a Leibniz graph with one Jacobiator]
Let $\Gamma$~be a Leibniz graph with one Jacobiator vertex $\Jac$.
From \eqref{EqJacFig} we see that expansion of $\Jac$ into a sum of three Kontsevich graphs means adding one new edge $w \to v$ (namely joining the internal vertices $w$ and $v$ within the Jacobiator).
Now, from $\Gamma$ construct three Kontsevich graphs by expanding $\Jac$ using \eqref{EqJacFig} and letting the edges which fall on $\Jac$ in $\Gamma$ be directed only to $v$ in every new graph.
Next, for each Kontsevich graph find the relabelling $\tau$ which brings it to its normal form and re-express the edge $w \to v$ using $\tau$.
Finally, out of the three normal forms of three graphs pick the minimal one.
By definition, the \emph{normal form} of the Leibniz graph $\Gamma$ is the pair: normal form of Kontsevich graph, that edge $\tau(w) \to \tau(v)$.
\end{define}

We say that a sum of Leibniz graphs is a \emph{skew Leibniz graph}~$\Alt\Gamma$ 
if it is produced from a given Leibniz graph~$\Gamma$
by alternation using formula~\eqref{EqMakeSkew}.

\begin{define}[Normal form of a skew Leibniz graph with one Jacobiator]\label{DefSkewNormalForm}
Likewise, the normal form of a skew Leibniz graph~$\Alt\Gamma$ is the minimum of the normal forms of Leibniz graphs (specifically, of the graph but not edge encodings) which are obtained from~$\Gamma$ by running over the  group of permutations of the sinks content.
\end{define}

\begin{lemma}[\cite{sqs15}]
In order to show that a sum~$S$ of weighted skew\/-\/sym\-met\-ric 
Kontsevich graphs vanishes for all Poisson structures $\cP$, 
it suffices 
to express~$S$ as a sum of skew Leibniz graphs: 
$
S = \Diamond\bigl(\cP,\Jac(\cP)\bigr)$.
\end{lemma}

\subsection{Formulation of the problem}
Let $\cP \mapsto \cP + \varepsilon \cQ (\cP) + \bar{o}(\varepsilon)$ be a deformation of bi-vectors that preserves their property to be Poisson at least infinitesimally on all affine manifolds:
$
\lshad\cP + \varepsilon\cQ + \bar{o}(\varepsilon), \cP + \varepsilon\cQ + \bar{o}(\varepsilon)\rshad = \bar{o}(\veps)$.
Expanding and equating the first order terms, we obtain the equation
$
\lshad\cP,\cQ (\cP)\rshad \doteq 0$
 {via} 
$ \lshad\cP ,\cP\rshad = 0$.
The language of Kontsevich graphs 
allows one to convert this infinite analytic problem within a given set\/-\/up $\bigl(\mathcal{N}^n,\mathcal{P}\bigr)$ in dimension~$n$ into a set of finite combinatorial problems whose solutions are universal for all Poisson geometries in all dimensions~$n<\infty$.

Our first task in this paper is 
to find the space of flows $\dot{\cP} = \cQ (\cP)$ which are encoded by the Kontsevich graphs on a fixed number of internal vertices~$k$, for $1\leqslant k\leqslant 4$. Specifically, we solve the graph equation
\begin{equation}\label{EqDiamond}
\lshad\cP,\cQ (\cP)\rshad = \Diamond\bigl(\cP,\Jac(\cP)\bigr)
\end{equation}
for the Kontsevich bi-vector graphs~$\cQ(\cP)$ and 
Leibniz graphs~$\Diamond$. 
We then 
factor out the Poisson\/-\/trivial and improper solutions, that is, we quotient out all bi\/-\/vector graphs that can be written in the form 
$
\cQ(\cP) = \lshad\cP ,X\rshad + \nabla\bigl(\cP,\Jac(\cP)\bigr)$,
where $X$~is a Kontsevich one\/-\/vector graph and $\nabla$~is a Leibniz bi\/-\/vector graph. (The bi\/-\/vectors $\lshad\cP, X\rshad$ make $\lshad\cP , \cQ(\cP)\rshad$ vanish since $\lshad\cP,\cdot\rshad$~is a differential. 
The improper graphs $\nabla (\cP,\Jac (\cP))$ 
vanish identically at all Poisson bi-vectors~$\cP$ on every affine manifold.

Before solving 
factorisation problem~\eqref{EqDiamond} 
with respect to the operator~$\Diamond$,
we must generate --\,e.g., iteratively as described below\,--
an ansatz for expansion of the right\/-\/hand side 
using skew Leibniz graphs with un\-de\-ter\-mi\-ned coefficients. 


\subsection{How to generate Leibniz graphs iteratively}\label{SecGenLeibnizIteratively}
The first step is to construct a 
layer of skew Leibniz graphs, 
that is, all skew Leibniz graphs which produce at least one graph in the input
(in the course of expansion of skew Leibniz graphs using formula~\eqref{EqMakeSkew} and then
in the course of expansion of every Leibniz graph at hand to a sum of Kontsevich graphs).
For a given Kontsevich graph in the input~$S$, one such Leibniz graph can be constructed 
by contracting an edge between two internal vertices 
so that the new vertex with three outgoing edges becomes 
the Jacobiator vertex. 
Note that these Leibniz graphs, which are designed to reproduce~$S$,
may also produce extra Kontsevich graphs that are not given in the input.
Clearly, if the set of Kontsevich graphs in~
$S$ coincides with the set of such graphs obtained by expansion of all the Leibniz graphs
in the ansatz at hand, then we are done: the extra graphs, not present in~$S$, are known to all cancel.
Yet it could very well be that it is not possible to express $S$ using only the Leibniz graphs 
from the set accumulated so far. 
Then we 
proceed by constructing the next 
layer of skew Leibniz graphs that reproduce at least one of the extra Kontsevich graphs 
(which were not present in~$S$ but which are produced by the graphs in the previously constructed layer(s) 
of Leibniz graphs).
In this way we 
proceed iteratively until no new Leibniz graphs are found;
of course, the overall number of skew Leibniz graphs on a fixed number of internal vertices and sinks 
is bounded from above so that the algorithm always terminates.
Note that the Leibniz graphs 
obtained in this way are the only ones that can in principle be involved in the vanishing mechanism for~$S$.



\begin{notation}\label{pNotationNHT}
Let $v$~be a graph 
vertex. Denote by~$N(v)$ the set of neighbours of~$v$, 
by~$H(v)$ the (possibly 
empty) set of arrowheads 
of oriented edges issued from the vertex~$v$, and
by~$T(v)$ the (possibly 
empty) set of tails 
for oriented edges pointing at~$v$.
For example, $\#N(\bullet)=2$, $\#H(\bullet)=2$, and $T(\bullet)=\varnothing$ 
for the top~$\bullet$ of the wedge graph~${\gets}\bullet{\to}$.
\end{notation}

\smallskip\noindent%
{\textbf{Algorithm}}
Consider a skew\/-\/symmetric sum~$S_0$ of oriented Kontsevich graphs
with real coefficients.
Let $S_{\text{total}} \mathrel{{:}{=}} S_0$ and create an empty table~$L$. 
We now describe the $i$th iteration of the algorithm~($i\geqslant 1$).\\[3pt]
\mbox{ }\quad{\textbf{\textsf{Loop~$\circlearrowright$}}}
Run over all Kontsevich graphs~$\Gamma$ in $S_{i-1}$: for each internal vertex~$v$ 
in a graph~$\Gamma$, run over all vertices~$w\in T(v)$ 
in the set of tails of oriented edges pointing at~$v$
such that $v\notin T(w)$ and $H(v)\cap H(w)=\varnothing$
for the sets of targets of oriented edges issued from~$v$ and~$w$.
Replace the edge $w\to v$ connecting $w$ to~$v$ by 
Jacobiator~\eqref{EqJacFig},
that is, by a single vertex~$\Jac$
with three outgoing edges
and such that $T(\Jac)= \bigl(T(v)\setminus w\bigr)\cup T(w)$ 
and $H(\Jac)=H(v)\cup (H(w)\setminus v) \mathrel{{=}{:}} \{a,b,c\}$.

Because we shall always expand the skew Leibniz graphs in what follows, 
we do not actually contract the edge $w\to v$ (to obtain a 
Leibniz graph explicitly) in this algorithm
but instead we continue working with the original Kontsevich graphs containing the distinct vertices~$v$ and~$w$.

For every edge that points at~$w$, redirect it to~$v$.
Sum over the three cyclic permutations that provide 
three possible ways to attach the three outgoing edges for~$v$ and~$w$ (excluding~$w\to v$) 
--\,now 
seen as the outgoing edges of the Jacobiator\,-- to the target vertices~$a$, $b$, and~$c$
depending on~$w$ and~$v$.
Skew\/-\/symmetrise\footnote{%
This algorithm can be modified so that it works for an input which is not skew, namely,
by replacing skew Leibniz graphs by ordinary Leibniz graphs (that is, by skipping the skew\/-\/symmetrisation).
For example, this strategy has been used in~\cite{sqs15,cpp} to show that the Kontsevich star 
product~$\star\mod\bar{o}(\hbar^4)$ is associative: although
the associator $(f\star g)\star h - f\star(g\star h) = \Diamond\bigl(\cP,\Jac(\cP)\bigr)$ 
is not skew, it does vanish for every Poisson structure~$\cP$.%
}
each of these three graphs with respect to the sinks
by applying for\-mu\-la~\eqref{EqMakeSkew}.


For every marked edge $w\to v$ 
indicating the internal 
edge in the Jacobiator vertex in a graph, replace each 
sum of the Kontsevich graphs which is skew with respect to the sink content by using the normal form of the respective skew Leibniz graph, see Definition~\ref{DefSkewNormalForm}.
%
%
If this skew Leibniz graph is not contained in~$L$, 
apply the Leibniz rule(s) for all the derivations acting on the Jacobiator vertex~$\Jac$. 
Otherwise speaking, 
sum over all possible ways to attach the incoming edges of the target~$v$ 
in the marked edge~$w\to v$ to its source~$w$ and target.
To each Kontsevich graph resulting from a skew Leibniz graph at hand assign 
the same undetermined coefficient, and add all these weighted Kontsevich graphs to the sum~$S_{i}$. 
Further, add a row to the table~$L$, that new row containing the normal form of
this skew Leibniz graph (with its coefficient
that has been made common to the Kontsevich graphs). 

By now, the new sum of Kontsevich graphs~$S_i$ is fully composed.
Having thus finished 
the current iteration over all graphs~$\Gamma$ in the set~$S_{i-1}$,
redefine the algebraic sum of weighted graphs~$S_{\text{total}}$ by subtracting from it the newly formed sum~$S_i$.
Collect similar terms in~$S_{\text{total}}$ and reduce this sum of Kontsevich graphs
modulo their skew\/-\/symmetry under swaps $L\rightleftarrows R$ of the edge ordering in every internal vertex,
so that all zero graphs (see Example~\ref{ExZeroGraph}) also get eliminated.
\hfill $\circlearrowright$~\textbf{\textsf{end loop}}

Increment~$i$ by~$1$ and repeat the iteration until the set of weighted (and skew) Leibniz graphs~$L$ stabilizes.
Finally, solve --\,with respect to the coefficients of skew Leibniz graphs\,--
the linear algebraic system obtained from the graph equation 
$S_{\text{total}} = 0$ for the sum of Kontsevich graphs which has been produced from its initial value~$S_0$ by running the iterations of the above algorithm.

\begin{example}
For the skew sum of Kontsevich graphs in the right-hand side of \eqref{EqJacFig}, the algorithm would produce just one skew Leibniz graph: namely, the Jacobiator itself.
\end{example}

\begin{example}[The $3$-\/wheel]
For the Kontsevich tetrahedral flow $\dot\cP=\cQ_{1:6/2}(\cP)$ on the spaces of Poisson bi\/-\/vectors~$\cP$, see~\cite{Ascona96,Kontsevich2017Bourbaki} and~\cite{tetra16,f16}, building a sufficient set of skew Leibniz graphs in the r.\/-\/h.\,s.\ of factorisation problem \eqref{EqDiamond} 
requires two iterations of the above algorithm: 11
Leibniz graphs are produced at the first step and 50~
more are added by the second step, making~61 in total.
One of the two known solutions of this factorisation problem~\cite{f16} then consists of 8~
skew 
Leibniz graphs (expanding to 27~Leibniz graphs). In turn, as soon as all the Leibniz rules acting on the Jacobiators are processed and every Jacobiator vertex is expanded via~\eqref{EqJacFig}, the right\/-\/hand side~$\Diamond\bigl(\cP,\Jac(\cP)\bigr)$ equals the sum of 39~Kontsevich graphs which are assembled into the 9~totally skew\/-\/symmetric terms in the left\/-\/hand side~$\lshad\cP,\cQ_{1:6/2}\rshad$. 
\end{example}

\begin{example}[The $5$-\/wheel]
Consider the factorisation problem $\lshad\cP, 
\Ori(\boldsymbol{\gamma}_5) \rshad = \Diamond\bigl(\cP,\Jac(\cP)\bigr)$ for the pentagon\/-\/wheel 
deformation $\dot\cP=\Ori(\boldsymbol{\gamma}_5) 
(\cP)$ of Poisson bi\/-\/vectors~$\cP$, see~\cite{Kontsevich2017private,Willwacher2017private} and~\cite{WeFactorize5Wheel}.
The ninety skew Kontsevich graphs encoding the bi\/-\/vector~$\Ori(\boldsymbol{\gamma}_5)
$ are obtained by taking all the admissible orientations of two $(5+1)$-\/vertex graphs~$\boldsymbol{\gamma}_5$, one of which is the pentagon wheel with five spokes, the other graph complementing the former to a cocycle in the non\/-\/oriented graph complex.
By running the iterations of the above algorithm for self\/-\/expanding construction of the 
Leibniz tri\/-\/vector graphs in this factorisation problem, we achieve 
stabilisation of the number of such graphs after the seventh iteration, see Table~\ref{Tab5WheelLeibniz} below.%
\begin{table}[htb]
\caption{The number of skew Leibniz graphs produced 
iteratively for~$\lshad\cP,
 \Ori(\boldsymbol{\gamma}_5 )   
\rshad$.}\label{Tab5WheelLeibniz}
\centerline{\begin{tabular}{l c c c c c c c c }
\br
No.\ iteration~$i$\strut & 1 & 2 & 3 & 4 & 5 & 6 & 7 & 8 \\
\mr
$\#$ of graphs\strut & 1518 & 14846 & 41031 & 54188 & 56318 & 56503 & 56509 & 56509 \\
\mr
\mbox{ } of them new\strut & all & +13328 & +26185 & +13157 & +2130 & +185 & +6 & none \\
\br
\end{tabular}
}
\end{table}
\end{example}


\section{Generating the Kontsevich multi\/-\/vector graphs}
\noindent%
Let us return to problem \eqref{EqDiamond}: it is the ansatz for bi\/-\/vector Kontsevich graphs $\cQ(\cP)$ with $k$ internal vertices, as well as the Kontsevich 1-vectors $\cX$ with $k-1$ internal vertices (to detect trivial terms $\cQ(\cP) = \schouten{\cP, \cX}$) which must be generated at a given $k$.
(At $1 \leqslant k \leqslant 4$, one can still expand with respect to \emph{all} the Leibniz graphs in the r.-h.s. of \eqref{EqDiamond}, not employing the iterative algorithm from \S\ref{SecGenLeibnizIteratively}.
So, a generator of the Kontsevich (and Leibniz) tri-vectors will also be described presently.)

The Kontsevich graphs corresponding to $n$-vectors are those graphs with $n$ sinks
(each containing the respective argument of $n$-\/vector) 
in which 
exactly one arrow comes into each sink, 
so that the order of the differential operator encoded by an $n$-\/vector graph 
equals one w.r.t.\ each argument, and which are totally skew-symmetric in their
$n$~arguments. 
Let us explain how one can economically obtain the set of one\/-\/vectors and
skew\/-\/symmetric bi\/- and tri\/-\/vectors with $k$~internal
vertices in three steps (including graphs with eyes $\bullet \rightleftarrows \bullet$
but excluding graphs with 
tadpoles). This approach can easily be extended to the construction of 
$n$-vectors with any~$n\geqslant1$. 

\subsection{\textbf{One\/-\/vectors}}
Each one\/-\/vector under study is encoded by a Kontsevich graph with one sink.
Since the sink has one incoming arrow, there is an 
internal vertex 
as the tail of this incoming arrow. 
The target of another edge  issued from this internal vertex can be any internal vertex other then itself. 

\smallskip
\noindent{\textit{Step~1.}} 
Generate all Kontsevich graphs on $k-1$ internal vertices and one sink
(i.e.\ graphs including those with eyes yet excluding those with tadpoles, 
and not necessarily of differential order one with respect to the sink content). 

\noindent%
{\textit{Step~2.}} For every such graph with $k-1$~internal vertices,
add the new sink and make it a target of the old sink, which itself becomes the
$k$th internal vertex.
Now run over the $k-1$ internal vertices excluding the old sink and --\,via the Leibniz rule\,-- make every such internal vertex the second target of the old sink.





\subsection{\textbf{Bi\/-\/vectors}}
There are two cases in the construction of bi\/-\/vectors encoded by the
Kontsevich graphs. At all~$k\geqslant1$ the first variant is referred to
those graphs with an internal vertex that has both sinks as targets.

\smallskip\noindent%
{\textit{Variant 1}: \textit{Step 1.}} Generate all $k$-vertex graphs on
$k-1$~internal vertices and one sink.

\noindent%
{\textit{Variant 1}: \textit{Step 2.}} For every such graph, 
add two new sinks and proclaim them as targets of the old sink.

Note that the obtained graphs are skew\/-\/symmetric. 

\smallskip
The second variant produces those graphs which contain two internal vertices such that 
one has the first sink as target and the other has the second sink as target.
The second target of either such internal vertex can be any internal vertex other then itself. Note that for $k = 1$ only 
the first variant applies.

\smallskip\noindent%
{\textit{Variant 2}: \textit{Step 1.}} Generate all $k$-vertex Kontsevich graphs on
$k-2$~internal vertices  and two sinks.
These sinks now become
the $(k-1)$th and $k$th~internal vertices.

\noindent{\textit{Variant 2}: \textit{Step 2.}} For every such graph, add two new sinks, make the first new sink a target of the first old sink and make the second 
new sink a target of the second old sink. Now run over the $k-1$ internal 
vertices excluding the first old sink, each time proclaiming an internal vertex 
the second target of the first old sink. 
Simultaneously, run over the $k-1$ internal vertices excluding the second old
sink and likewise, declare an internal vertex to be the second target of the second 
old sink.

\noindent{\textit{Variant 2}: \textit{Step 3.}} 
Skew\/-\/symmetrise each graph with respect to the content of two sinks using~\eqref{EqMakeSkew}.


\subsection{\textbf{Tri\/-\/vectors}}
For $k\geqslant 3$, there exist two variants of tri\/-\/vectors. The first variant
at all~$k\geqslant2$ yields those Kontsevich graphs with two internal vertices such that one has 
two of the three sinks as its targets 
while another internal vertex has the third sink as one of its targets. 
The second target of this last vertex can be any internal vertex other then itself.
The second variant contains those graphs with three internal vertices such that
the first one has the first sink as a target, the second one has the second 
sink as a target, and the third one has the third sink as a target. 
For each of these three internal vertices with a sink as target, 
the second target can be any internal vertex other then itself.

\smallskip\noindent%
{\textit{Variant 1}: \textit{Step 1.}} Generate all 
$k$-vertex Kontsevich graphs on $k-2$~internal vertices and two sinks.

\noindent{\textit{Variant 1}: \textit{Step 2.}} For every such graph,
add three new sinks, make the first two new sinks the targets of the 
first old sink and make the third new sink a target of the second old sink.
Now run over the $k-1$ internal vertices excluding the second old sink and 
every time declare an internal vertex the second target of the second old sink.

\noindent{\textit{Variant 1}: \textit{Step 3.}}
Skew\/-\/symmetrise all 
graphs at hand by applying formula~\eqref{EqMakeSkew} to each of~them.

Note that for $k=1$ there are no tri\/-\/vectors encoded by Kontsevich graphs 
and also note that for $k=2$ 
only 
the first variant applies.

\smallskip\noindent%
{\textit{Variant 2}: \textit{Step 1.}} Generate all 
Kontsevich graphs on $k-3$~internal vertices and three sinks.

\noindent
{\textit{Variant 2}: \textit{Step 2.}} For every such graph, add three new sinks, make the first new sink a target of the first old sink, 
make the second new sink a target of the second old sink and make the third 
new sink a target of the third old sink. Now run over the $k-1$~internal
 vertices excluding the first old sink and declare 
every such internal vertex 
the second target of the first old sink. Independently, run over the $k-1$~internal 
vertices excluding the second old sink and declare each internal vertex to be 
the second target of the second old sink. Likewise, run over the $k-1$ internal 
vertices excluding the third old sink and declare each internal vertex to 
be the second target of the third old sink.

\noindent{\textit{Variant 2}: \textit{Step 3}}
Skew\/-\/symmetrise all the graphs at hand using~\eqref{EqMakeSkew}.


\subsection{\textbf{Non\/-\/iterative generator of the Leibniz $n$\/-\/vector graphs}}
The following algorithm generates all Leibniz graphs with a prescribed number of internal vertices and sinks. 
{Note that not only multi\/-\/vectors, but also all graphs of arbitrary differential order with respect to the sinks can be generated this way.}

\smallskip\noindent%
{\textit{Step 1:}} Generate all Kontsevich graphs of prescribed type on $k-1$ internal vertices and $n$~sinks, e.g., all $n$\/-\/vectors.

\noindent%
{\textit{Step 2:}} Run through the set of these Kontsevich graphs and in each of them, run through the set of its internal vertices~$v$.
For every vertex~$v$ 
do the following: re\/-\/enumerate the internal vertices so that this vertex is enumerated by~$k-1$.
This 
vertex already targets two 
vertices, $i$ and~$j$, where $i<j<k-1$.
Proclaim the last, ($k-1$)th vertex to be the placeholder of the Jacobiator 
(see~\eqref{EqJacFig}), so we must still add the third arrow.
Let a new index~$\ell$ run up to~$i-1$ starting at~$n$ if 
only the $n$\/-\/vectors are produced.\footnote{%
If we want to generate not only $n$\/-\/vectors but all graphs of arbitrary differential orders, then we let $\ell$~start at~$0$ (so that the sinks are included).}
For every admissible value of~$\ell$, generate a new graph where the $\ell$th~vertex is proclaimed the third target of the Jacobiator vertex~$k-1$.
(
Restricting~$\ell$ by~$<i$, we reduce the number of possible repetitions in the set of Leibniz graphs. Indeed, for every triple $\ell<i<j$, the same Leibniz graph in which the Jacobiator stands on those three vertices would be produced from the three Kontsevich graphs: namely, those in which the $(k-1)$th~vertex targets at the~$\ell$th and~$i$th, at the~$\ell$th and~$j$th, and at the~$i$th and~$j$th vertices. 
In these three cases it is the~$j$th, $i$th, and $\ell$th~vertex, respectively, which would be appointed by the algorithm as the Jacobiator's third target.)


We use this algorithm to generate the Leibniz tri\/-{} and bi\/-\/vector graphs:
to establish Theorem~\ref{ThMainFew}, 
we list all possible terms in the right\/-\/hand side of 
factorisation problem~\eqref{EqDiamond} at~$k\leqslant 4$ and
then we filter out the improper bi\/-\/vectors in 
the found 
solutions~$\cQ(\cP)$.

\begin{rem}
There are 
at least 265,495
Leibniz 
graphs on $3$~sinks and $6$~internal vertices of which one is the Jacobiator vertex. When compared with Table~\ref{Tab5WheelLeibniz} on p.~\pageref{Tab5WheelLeibniz}, this estimate 
suggests why at large $k\gtrsim 5$, the 
breadth\/-\/first\/-\/search iterative algorithm from~
\S\ref{SecGenLeibnizIteratively} generates a smaller number of the Leibniz  tri\/-\/vector
graphs, namely, only the ones which can in principle be involved in the factorisation under study.
\end{rem}

\section{Main result} 

\begin{theor}[$k\leqslant4$]\label{ThMainFew}
The few-vertex solutions of problem \eqref{EqDiamond} are these \textup{(}note that disconnected Kontsevich graphs in $\cQ(\cP)$ are allowed\,\textup{!):}
\mbox{ %
}\\[1pt]
$\bullet$\ $k=1$\textup{\textbf{:}} The dilation $\dot{\cP} = \cP$ is a 
unique, nontrivial and proper solution.\\
$\bullet$\ $k=2$\textup{\textbf{:}} No solutions exist \textup{(}in particular, neither trivial nor improper\textup{)}. \\
$\bullet$\ $k=3$\textup{\textbf{:}} 
There are no solutions\textup{:} 
neither Poisson\/-\/trivial nor Leibniz bi\/-\/vectors.\\ 
$\bullet$\ $k=4$\textup{\textbf{:}} A 
unique nontrivial and proper solution is the Kontsevich tetrahedral 
flow~$\cQ_{1:\frac{6}{2}}(\cP)$ from Example~\textup{\ref{Ex3Wheel}}
\textup{(}see~\textup{\cite{Ascona96,Kontsevich2017Bourbaki}} and~\textup{\cite{tetra16,f16})}.
There is a one\/-\/dimensional space of Poisson trivial \textup{(}still proper\textup{)} solutions $\lshad\cP,X\rshad$\textup{;} 
the Kontsevich $1$-\/vector~$X$ on three internal vertices is drawn in~\textup{\cite[App.\:F]{f16}}. 
Intersecting with the former by $\{0\}$, there is a three\/-\/dimensional space of improper \textup{(}still Poisson\/-\/nontrivial\textup{)} solutions of the form~$\nabla\bigl(\cP,\Jac(\cP )\bigr)$. 

None of the solutions~$\cQ(\cP)$ known so far contains any $2$\/-\/cycles \textup{(}or ``eyes"~$\bullet\rightleftarrows\bullet$\textup{)}.\footnote{Finding 
solutions~$\cQ(\cP)$ with tadpoles or extra sinks --\,with fixed arguments\,-- 
is a separate 
problem
.} 
\end{theor}


We now 
report a 
classification of Poisson bi\/-\/vector symmetries $\dot\cP=\cQ(\cP)$ which are 
given by those Kontsevich graphs~$\cQ=\Ori(\gamma)$ on $k$~internal vertices  that can be obtained at $5\leqslant k\leqslant 9$ by orienting $k$-\/vertex connected graphs~$\gamma$ without double edges. By construction, this extra assumption keeps 
only those 
Kontsevich graphs which may not contain~eyes.

We first find 
such graphs~$\gamma$ that satisfy~$\Id(\gamma)=0$, then we exclude the coboundaries $\gamma = \Id(\gamma')$ for some graphs~$\gamma'$ on~$k-1$ vertices and $2k-3$~edges.


\begin{theor}[$5\leqslant k\leqslant 8$]\label{ThMainMany}
Consider the vector space of non\/-\/oriented connected graphs on $k$~vertices and $2k-2$~edges,
without tadpoles and without multiple edges. 
All nontrivial $\Id$-\/cocycles for $5\leqslant k\leqslant 8$ are exhausted by the following ones\textup{:}
\\[1pt]
$\bullet$\ $k=5,7$\textup{\textbf{:}} No solutions. \\ 
$\bullet$\ $k=6$\textup{\textbf{:}}\label{Where5WheelIs} %
{\unitlength=1mm
\begin{picture}(0,0)(-15,1.5)
\put(65,0){$\boldsymbol{\gamma}_5 = {}$}
\put(85,0){
{\unitlength=0.3mm
\begin{picture}(55,53)(5,-5)
\put(27.5,8.5){\circle*{3}}
\put(0,29.5){\circle*{3}}
\put(-27.5,8.5){\circle*{3}}
\put(-17.5,-23.75){\circle*{3}}
\put(17.5,-23.75){\circle*{3}}
\qbezier
(27.5,8.5)(0,29.5)(0,29.5)
\qbezier
(0,29.5)(-27.5,8.5)(-27.5,8.5)
\qbezier
(-27.5,8.5)(-17.5,-23.75)(-17.5,-23.75)
\qbezier
(-17.5,-23.75)(17.5,-23.75)(17.5,-23.75)
\qbezier
(17.5,-23.75)(27.5,8.5)(27.5,8.5)
\put(0,0){\circle*{3}}
\qbezier
(27.5,8.5)(0,0)(0,0)
\qbezier
(0,29.5)(0,0)(0,0)
\qbezier
(-27.5,8.5)(0,0)(0,0)
\qbezier
(-17.5,-23.75)(0,0)(0,0)
\qbezier
(17.5,-23.75)(0,0)(0,0)
\end{picture}
}
}
\put(95,0){${}+\dfrac{5}{2}$}
\put(114,0){
{\unitlength=0.4mm
\begin{picture}(50,30)(0,-4)
\put(12,0){\circle*{2.5}}
\put(-12,0){\circle*{2.5}}
\put(25,15){\circle*{2.5}}
\put(-25,15){\circle*{2.5}}
\put(-25,-15){\circle*{2.5}}
\put(25,-15){\circle*{2.5}}
\put(-12,0){\line(1,0){24}}
\put(-25,15){\line(1,0){50}}
\put(-25,-15){\line(1,0){50}}
\put(-25,-15){\line(0,1){32}}
\put(25,15){\line(0,-1){32}}
\qbezier
(25,15)(12,0)(12,0)
\qbezier
(-25,15)(-12,0)(-12,0)
\qbezier
(-25,-15)(-12,0)(-12,0)
\qbezier
(25,-15)(12,0)(12,0)
\put(-12.5,17){\oval(25,10)[t]}
\put(12.5,-17){\oval(25,10)[b]}
\put(0,2){\line(0,1){11}}
\put(0,-2){\line(0,-1){11}}
\end{picture}
}%
}
\end{picture}%
}%
A unique solution\footnote{There are only 12
admissible graphs to build cocycles from; of these 12, as many as 6 
are zero graphs. This count shows 
to what extent the number of graphs decreases if one restricts to only the flows~$\cQ=\Ori(\gamma)$ obtained from cocycles~$\gamma\in\ker(\Id)$ in the non\/-\/oriented graph complex.} 
is given by the Kontse-\\
\parbox[t]{86mm}{vich\/--\/Willwacher pentagon\/-\/wheel cocycle \textup{(}
see Example~\textup{\ref{Ex5Wheel}).} 
The established 
factorisation $\lshad\cP,\Ori(\boldsymbol{\gamma}_5)\rshad=\Diamond\bigl(\cP,\Jac(\cP)\bigr)$ 
$\lefteqn{\text{will be addressed in a separate paper \textup{(}see~\textup{\cite{WeFactorize5Wheel})}.}}$}%
\\[0.5pt]
$\bullet$\ $k=8$\textup{\textbf{:}} The only solution~$\boldsymbol{\gamma}_7$ consists of the heptagon\/-\/wheel and $45$~other graphs \textup{(}see Table~\textup{\ref{Tab1+45}}, in which the coefficient of heptagon graph is~$\mathbf{1}$ in bold,
and~\textup{\cite{JNMP2017})}.%
\begin{table}[htb]
\caption{The heptagon\/-\/wheel graph cocycle~$\boldsymbol{\gamma}_7$.}\label{Tab1+45}
\centerline{{\upshape\footnotesize
\begin{tabular}{l r l r} 
\br
Graph encoding & Coeff.\ \  & Graph encoding & Coeff.\\
\mr
16	17	18	23	25	28	34	38	46	48	57	58	68	78	& 	$\mathbf{1}\phantom{/1}$		& 	12	13	18	25	26	37	38	45	46	47	56	57	68	78	&	$-7\phantom{/1}$	\\
12	14	18	23	27	35	37	46	48	57	58	67	68	78	& 	$-21/8$	& 	12	14	16	23	25	36	37	45	48	57	58	67	68	78	&	$77/8$	\\
13	14	18	23	25	28	37	46	48	56	57	67	68	78	& 	$-77/4$	& 	13	16	17	24	25	26	35	37	45	48	58	67	68	78	&	$-7\phantom{/1}$	\\
12	13	15	24	27	35	36	46	48	57	58	67	68	78	& 	$-35/8$	& 	14	15	17	23	26	28	37	38	46	48	56	57	68	78	&	$49/4$	\\
12	13	18	24	26	37	38	46	47	56	57	58	68	78	& 	$49/8$	& 	12	16	18	27	28	34	36	38	46	47	56	57	58	78	&	$-147/8$	\\
14	17	18	23	25	26	35	37	46	48	56	58	67	78	& 	$77/8$	& 	12	15	16	27	28	35	36	38	45	46	47	57	68	78	&	$-21/8$	\\
12	13	18	26	27	35	38	45	46	47	56	57	68	78	& 	$-105/8$		& 	12	14	18	23	27	35	36	45	46	57	58	67	68	78	&	$-35/8$	\\
12	14	18	23	27	36	38	46	48	56	57	58	67	78	& 	$7/8$	& 	14	15	16	23	26	28	37	38	46	48	57	58	67	78	&	$-49/4$	\\
12	14	15	23	27	35	36	46	48	57	58	67	68	78	& 	$35/8$	& 	12	15	18	23	28	34	37	46	48	56	57	67	68	78	&	$105/8$	\\
12	13	14	27	28	36	38	46	47	56	57	58	68	78	& 	$-49/8$	& 	12	14	17	23	26	37	38	46	48	56	57	58	68	78	&	$-49/8$	\\
12	13	18	25	27	34	36	47	48	56	58	67	68	78	& 	$35/4$	& 	12	16	18	25	27	35	36	37	45	46	48	57	68	78	&	$49/1\lefteqn{6}$	\\
12	13	14	25	26	36	38	45	47	57	58	67	68	78	& 	$-119/1\lefteqn{6}$	& 	12	13	18	25	27	35	36	46	47	48	56	57	68	78	&	$7\phantom{/1}$	\\
12	13	15	24	28	36	38	47	48	56	57	67	68	78	& 	$49/8$	& 	12	14	18	25	28	34	36	38	47	57	58	67	68	78	&	$-7\phantom{/1}$	\\
12	13	14	23	28	37	46	48	56	57	58	67	68	78	& 	$77/4$	& 	12	16	18	25	27	35	36	37	45	46	48	58	67	78	&	$-77/1\lefteqn{6}$	\\
12	15	17	25	26	35	36	38	45	47	48	67	68	78	& 	$-49/8$	& 	12	14	18	23	27	35	38	46	47	57	58	67	68	78	&	$77/4$	\\
13	15	18	24	26	28	37	38	46	47	56	57	68	78	& 	$-49/4$	& 	12	14	15	23	27	36	38	46	48	57	58	67	68	78	&	$35/2$	\\
13	14	18	25	26	28	36	38	47	48	56	57	67	78	& 	$-49/4$	& 	12	13	18	25	27	34	36	46	48	57	58	67	68	78	&	$-105/8$	\\
12	14	18	23	28	35	37	46	48	56	57	67	68	78	& 	$-7\phantom{/1}$ &	12	15	16	25	27	35	36	38	46	47	48	57	68	78	&	$-7\phantom{/1}$	\\
12	14	18	23	28	36	38	46	47	56	57	58	67	78	& 	$-7\phantom{/1}$		& 	12	13	16	25	28	34	37	47	48	57	58	67	68	78	&	$-147/1\lefteqn{6}$	\\
12	15	16	25	27	35	36	38	46	47	48	58	67	78	& 	$49/8$	& 	12	13	17	25	26	35	37	45	46	48	58	67	68	78	&	$-77/4$	\\
12	14	18	23	28	36	37	46	47	56	57	58	68	78	& 	$49/8$		& 	12	14	17	23	27	35	38	46	48	57	58	67	68	78	&	$-49/8$	\\
12	13	15	26	27	35	36	45	47	48	58	67	68	78	& 	$-7\phantom{/1}$	& 	12	13	15	26	28	35	37	45	46	47	58	67	68	78	&	$-7/4$	\\
12	13	18	24	28	35	38	46	47	57	58	67	68	78	& 	$7\phantom{/1}$	& 	12	14	18	23	26	36	38	47	48	56	57	58	67	78	&	$-7\phantom{/1}$ 	\\
\br
\end{tabular}}}
\end{table}
\end{theor}

\begin{rem} 
The wheel graphs are built of 
triangles. The differential~$\Id$ cannot produce any triangle since 
multiple edges are not allowed. Therefore, all wheel cocycles are nontrivial. Note also that every wheel graph with $2\ell$~spokes is invariant under a mirror reflection with respect to a diagonal consisting of two edges attached to the centre. Hence there exists an edge permutation that swaps $2\ell-1$~pairs of edges. By footnote~\ref{FootZeroGraph} such graph equals~zero.
\end{rem}

\appendix
\section{How the orientation mapping~$\Ori$ is calculated%
}\label{SecOrient}
\noindent%
The algorithm lists all ways in which a given non-oriented graph can be oriented in such a way that it becomes a Kontsevich graph on two sinks. It consists of two steps: 
\begin{enumerate}
	\item choosing 
	the source(s) of the two arrows pointing at the first and second sink, respectively;
	\item orienting the edges between the internal vertices in all admissible ways, so that only Kontsevich graphs are obtained. 
\end{enumerate}

\smallskip\noindent%
{\textit{Step 1.}} Enumerate the $k$ vertices of a given non-oriented,
connected graph using~$2$, $\dots$, $k+1$. They become the internal vertices 
of the oriented graph. Now add the two sinks to the non\/-\/oriented graph, the sinks
enumerated using~$0$ and~$1$. Let~$a$ and~$b$ be a non\/-\/strictly ordered 
($a\leqslant b$) pair of internal vertices in the graph. Extend the graph by oriented edges $a\rightarrow 0$ and $b\rightarrow 1$ from vertices~$a$ and~$b$ to the sinks~$0$ and~$1$, respectively. 

\begin{rem} The choice of such a base pair, that is, the vertex or vertices from which two arrows are issued to the sinks, is an external input in the orientation procedure. 
Let us agree that if, at any step of the algorithm, a contradiction is achieved so that a graph at hand cannot be of Kontsevich type, the oriented graph draft is discarded; one proceeds with the next options in that loop, or if the former loop is finished, with the next level-up loops, or -- having returned to the choice of base vertices -- with the next base. In other words, we do not exclude in principle a possibility to have no admissible orientations for a particular choice of the base for a given non-oriented graph.
\end{rem}

\begin{notation}
Let $v$ be an internal vertex. 
Recalling from p.~\pageref{pNotationNHT} the notation for the set~$N(v)$ of neighbours of~$v$, the (initially empty) set~$H(v)$ of arrowheads of oriented edges issued from the vertex~$v$, and the (initially empty) set~$T(v)$ of tails for oriented edges pointing at $v$, we now put by definition
$F(v) \mathrel{{:}{=}} N(v)\setminus (H(v)\cup T(v))$. In other words, $F(v)$~is the subset of neighbours connected with~$v$ by a non\/-\/oriented edge.
\end{notation}

\noindent
{\textit{Step~2.1. Inambiguous orientation of \textup{(}some\textup{)} edges.}}
Here we use that every internal vertex of a Kontsevich graph should be the tail of 
exactly two outgoing arrows. 
We run 
over the set of all internal vertices~$v$. For every vertex such that 
the number of elements $\# H(v)=2$
, proclaim $T(v) \mathrel{{:}{=}} N(v)\setminus H(v)$, whence $F(v)=\varnothing$.
If for a vertex~$v$ we have that $\# H(v)=1$ and $\# F(v)=1$, then include $F(v)\hookrightarrow H(v)$, 
that is, convert 
a unique non\/-\/oriented edge touching~$v$ into an outgoing edge 
issued from 
this vertex. 
If $\# H(v)=0$ and $\# F(v)=2$, also include $F(v)\hookrightarrow H(v)$, 
effectively making both non\/-\/oriented edges outgoing from~$v$.

Repeat the three parts of Step~2.1 while any of the sets~$F(v)$, $T(v)$, or
~$S(v)$
is modified for at least one internal vertex~$v$ unless a contradiction is revealed. 
Summarising, Step~2.1 amounts to finding the edge orientations which are implied by the choice of the base pair $a$,~$b$ and by all the orientations of edges fixed earlier.

\noindent
{\textit{Step 2.2. Fixing the orientation of \textup{(}some\textup{)} remaining edges.}} 
Choose an internal ver\-tex~$v$ such that $H(v)<2$ and such that $H(v)\neq\varnothing$ or $T(v)\neq\varnothing$, that is, choose a vertex that is not yet equipped with two outgoing edges and that is attached to an oriented edge. 
If~$\# H(v)=1$, then run over the non\/-\/empty set~$F(v)$: for every vertex~$w$ in~$F(v)$, 
include $\{w\}\hookrightarrow H(v)$ and start over at Step~2.1.
Otherwise, i.e.\ if $H(v)=\varnothing$, run over all ordered pairs $(u,v)$ of vertices in the set~$F(v)$: for every such pair, make $H(v) \mathrel{{:}{=}} \{u,w\}$ and start over at Step~2.1.

By realising Steps~1 and~2 we accumulate the sum of fully oriented Kontsevich graphs.

\ack
A.\,V.\,Kiselev thanks the Organising committee of the international 
conference ISQS'25 on integrable systems and quantum symmetries 
(6--10~June 2017 in $\smash{\text{\v{C}VUT}}$ Prague, Czech Republic) for a warm atmosphere 
during the meeting. 
The authors are grateful 
to 
M.~Kontsevich and T.~Willwacher for helpful discussion.
We also thank 
Center for Information Technology of the University of Groningen 
for providing access to 
\textsf{Peregrine} high performance computing cluster.\quad
This research 
was supported in part by 
JBI~RUG project~106552 (Groningen, The Netherlands).
A part of this research was done while R.\,Buring and A.\,V.\,Kiselev were 
visiting at the~IH\'ES (Bures\/-\/sur\/-\/Yvette, France) and A.\,V.\,Kiselev 
was visiting at the~MPIM (Bonn, Germany).

\end{document}